\documentstyle[12pt,amsfonts]{article}

\newtheorem{theorem}{Theorem}[section]
\newtheorem{lemma}[theorem]{Lemma}
\newtheorem{prop}[theorem]{Proposition}

\newcommand{\Z}{{\mathbb{Z}}}
\newcommand{\R}{{\mathbb{R}}}

\newcommand{\V}{{\mathbb{F}}_{2}^{d}}

\newcommand{\pol}{F_{2}[u_{1},\ldots ,u_{d}]}
\newcommand{\lar}{{\cal R}_{d}^{(2)}} 

\title{Isomorphism rigidity of algebraic ${\mathbb{Z}}^{d}$-actions}
\author{Siddhartha Bhattacharya}
\begin{document}
\maketitle
\begin{abstract}
An algebraic ${\Z}^{d}$-action is an action of 
 ${\Z}^{d}$ on a compact abelian group 
$X$ by continuous automorphisms
of $X$. We prove that that for $d\ge 8$, there exist mixing
zero entropy algebraic ${\Z}^{d}$-actions which do not 
exhibit isomorphism rigidity property. 
\end{abstract}
\section{Introduction}
An algebraic ${\Z}^{d}$-action is an action 
$\alpha : {\bf n}\rightarrow \alpha({\bf n})$ of 
 ${\Z}^{d}$ on a compact abelian group $X$ by continuous automorphisms
of $X$. It is easy to see that any such action preserves $\lambda_{X}$, 
the Haar measure on $X$. If $\alpha$ is a homomorphism from
$ {\Z}^{d}$ to $GL(n, {\Z})$ for some $n\ge 1$, then the natural
action of $\alpha({\Z}^{d})$ on ${\R}^{n}$ induces an 
algebraic ${\Z}^{d}$-action on 
${\mathbb{T}}^{n}\cong {\mathbb{R}}^{n}/{\mathbb{Z}}^{n}$. Another
class of examples is given by group shifts : let $F$ be a finite abelian
group and let $S$ be the shift action of ${\Z}^{d}$ on
$F^{{\Z}^{d}}$ defined by
$$ S({\bf m})(x)({\bf n}) = x({\bf n} + {\bf m})\ \ \forall 
{\bf m},{\bf n}\in {\Z}^{d}.$$
A group shift is a closed shift invariant subgroup $X\subset F^{{\Z}^{d}}$,
together with the shift action of ${\Z}^{d}$ restricted to $X$.

An algebraic ${\Z}^{d}$-action $(X,\alpha)$ is said to be {\it irreducible}
if $X$ does not admit proper closed $\alpha$-invariant infinite 
subgroups. If $(X,\alpha)$ and $(Y,\beta)$ are two algebraic 
${\Z}^{d}$-actions
and $f : X\rightarrow Y$ is a measurable map then $f$ is said to
be a {\it measurable conjugacy\/} if $f$ is a measure space
isomorphism from $(X,\lambda_{X})$ to $(Y,\lambda_{Y})$ and 
for all ${\bf n} \in {\Z}^{d}$, 
$f\circ\alpha({\bf n}) = \beta({\bf n})\circ f$ a.e. $\lambda_{X}$.
An {\it algebraic conjugacy\/} from $(X,\alpha)$ to $(Y,\beta)$
is a continuous isomorphism $\theta$ from $X$ to $Y$, which satisfies 
 $\theta\circ\alpha({\bf n}) = \beta({\bf n})\circ \theta$ 
for all ${\bf n}$ in ${\Z}^{d}$. Two 
algebraic ${\Z}^{d}$-actions $(X,\alpha)$ and $(Y,\beta)$ are 
said to be {\it measurably conjugate \/} 
(resp. {\it algebraically conjugate \/}) if there is a 
measurable conjugacy (resp. algebraic conjugacy) from
$(X,\alpha)$ to $(Y,\beta)$. If $X$ and $Y$ are compact abelian groups
and $f : X \rightarrow Y$ is a measurable map, then $f$ is said to be
{\it affine \/} if there exists  an element $c\in Y$ and a 
continuous surjective group homomorphism $\theta : X\rightarrow Y$
such that $f(x) = c + \theta(x)$ a.e. $\lambda_{X}$.

\medskip
Recently, in \cite{KKS} and \cite{KS2} it has been shown  that the measurable 
orbit structure of a certain class of mixing zero entropy
algebraic ${\Z}^{d}$-actions exhibit strong rigidity properties. 
More specifically, it has been proved that if $(X,\alpha)$ and $(Y,\beta)$ 
are two  algebraic ${\Z}^{d}$-actions such that
the actions $\alpha$ and $\beta$ are irreducible, expansive and mixing
then every measurable conjugacy from $(X,\alpha)$ to $(Y,\beta)$ 
is an affine map (cf. \cite{KS2}, Corollary 1.2). 
 The question whether this form of 
rigidity occurs for all mixing zero entropy algebraic 
${\Z}^{d}$-actions has been raised by several authors.
In various degrees of generality,
several questions and conjectures about this aspect of 
mixing zero entropy algebraic 
${\Z}^{d}$-actions can be found in \cite{KS1}, 
\cite{KS2} and \cite{QT} .
All these questions can be viewed as special cases of the
following more general conjecture due to K. Schmidt 
(cf. \cite{Sc2}, Conjecture 3.5).

\medskip
\noindent
{\bf Conjecture.}
 {\it Let $d > 1$, and let $\alpha$ and $\beta$ be mixing 
algebraic ${\Z}^{d}$-actions on compact abelian groups $X$ and
$Y$, respectively. If $h(\alpha) = 0$, and if
$\phi : X\rightarrow Y$ is a measurable conjugacy of $\alpha$ and 
 $\beta$, then  $\phi$ is $\lambda_{X}$ a.e. equal to an affine map.
In particular, measurable conjugacy implies algebraic conjugacy.
\/}

\medskip
In this note we give a counter-example to the above conjecture.
More specifically we prove the following result.

\begin{theorem}\label{th}
For every $d\ge 8$, there exists a mixing zero entropy 
algebraic ${\Z}^{d}$-action $\alpha$ on a compact zero dimensional
abelian group $X$, and a non-affine homeomorphism $f : X\rightarrow X$,
such that $f$ preserves the Haar measure on $X$ and commutes with
the action $\alpha$.
\end{theorem}
\section{Markov subgroups}
For any $d\ge 1$, by $X_{d}$ we denote the group $({\Z/2\Z})^{{\Z}^{d}}$,
equipped with pointwise addition and the topology of pointwise 
convergence. It is easy to see that 
 $X_{d}$ is a compact zero dimensional abelian 
group. By $S$ we denote shift action of ${\Z}^{d}$ 
on $X_{d}$ defined by 
$$ S({\bf j})(x) ({\bf i}) = x({\bf i} +{\bf j})\ \ \forall 
{\bf i},{\bf j}\in {\mathbb{Z}}^{d}.$$
A {\it Markov subgroup\/} of $X_{d}$ is a closed subgroup
which is invariant under the shift action.
 In \cite{KS1} it was shown that
the dynamics of the shift action on  Markov subgroups can be 
studied using algebraic methods. We briefly recall the results that are
needed in our construction. For proofs, the reader is referred to 
\cite{Sc1}, Theorem 6.5 and Proposition 19.4.

\medskip  
Let ${\mathbb{F}}_{2}$ denote the field with two elements and for 
$d\ge 1$, let ${\cal R}_{d}^{(2)}$ denote the group-ring of
${\Z}^{d}$ with coefficients in ${\mathbb{F}}_{2}$. The ring 
${\cal R}_{d}^{(2)}$ can be identified with 
${\mathbb F}_{2}[u_{1}^{\pm 1},\ldots , u_{d}^{\pm 1}]$, 
the ring of Laurent
polynomials in $d$ commuting variables with coefficients in
${\mathbb F}_{2}$. Every element $p\in {\cal R}_{d}^{(2)}$
is written as
$$ p = \sum_{{\bf m}\in {\Z}^{d}} c_{p}({\bf m})u^{{\bf m}},$$
with $u^{{\bf m}} = u^{m_{1}}\cdots u^{m_{d}}$ and 
$c_{p}({\bf m}) \in {\mathbb F}_{2}$, where 
$c_{p}({\bf m}) = 0$ for all but finitely many ${\bf m}$.
For any $d\ge 1$, the group $X_{d}$ can be viewed as a 
${\cal R}_{d}^{(2)}$-module via the operation
$ p\cdot x = \sum c_{p}({\bf m})S({\bf m})(x).$
For any ideal $I\subset{\cal R}_{d}^{(2)}$ we define 
$X(I)\subset X_{d}$ by
$$X(I)
= \{ x\in X_{d}\ |\ p\cdot x = 0 \ \forall p\in I\}.$$
It is easy to see that $X(I)$ is a Markov subgroup 
of $X_{d}$. Conversely, given
any Markov subgroup $H$ of $X_{d}$, we define an ideal
$I(H)\subset {\cal R}_{d}^{(2)}$ by 
$$I(H) = \{ p\in {\cal R}_{d}^{(2)}\ |\ p\cdot x = 0 \ \ \forall
x\in H\}.$$
Using duality theory of compact abelian groups it can be shown that
for any ideal $J\subset {\cal R}_{d}^{(2)}$ and for any 
Markov subgroup $H\subset X_{d}$,
$I(X(J)) = J$ and $X(I(H)) = H$. Hence the correspondence
$H\mapsto I(H)$ is an order reversing bijection from
the set of all Markov subgroups of $X_{d}$ to 
the set of all ideals in ${\cal R}_{d}^{(2)}$. 
\begin{prop}\label{mar}
Let $d\ge 1$ and let $H\subset X_{d}$ be a Markov subgroup.
\begin{enumerate}
\item{The action $(H,S)$ has zero entropy if and only if
$H$ is a proper subgroup.}
\item{ If $I(H)$ is a prime ideal then the action $(H,S)$ is mixing 
if and only if for every non-zero ${\bf m}$, $u^{{\bf m}} - 1$ 
does not lie in  $I(H)$.}
\end{enumerate}
\end{prop}
\section{Binary linear codes}
 A {\it binary linear code \/} of length $d$ is a subspace $C$ of 
${\mathbb{F}}_{2}^{d}$. For any ${\bf v}\in \V$,
by $|{\bf v}|$ we denote the number of non-zero coordinates of ${\bf v}$.
A set $A\subset \V$  is said to be {\it even \/}(resp. {\it doubly even\/}) 
if for every element ${\bf v}$ of $A$, $|{\bf v}|$ is divisible by
2 (resp. divisible by 4). 
If ${\bf v}$ and ${\bf w}$ are two elements of
 $\V$ then
their dot product ${\bf v}\cdot {\bf w}$ is defined by 
${\bf v}\cdot {\bf w} = \sum v_{i}w_{i}$. For any set $A\subset \V$, by
$A^{\perp}$ we denote the binary linear code defined by
$$  A^{\perp} = \{ {\bf v}\in \V\ |\ {\bf v}\cdot {\bf w} = 
0 \ \ \forall {\bf w} \in A\}.$$
A set $A\subset \V$ is said to be {\it self orthogonal \/ } if 
$A \subset A^{\perp}$.

\medskip
{\bf Example 1 : } 
For any $d\ge 2$, let $E_{d}\subset \V$ be the
subspace consisting of all ${\bf v}$ such that 
$|{\bf v}|$ is even. If ${\bf 1}$ denotes the vector 
 $(1,1,\ldots ,1)$ in $\V$ then $E_{d}^{\perp}$ is an
one dimensional subspace of $\V$, consisting of 
${\bf 1}$ and ${\bf 0}$.

\medskip
{\bf Example 2 :}
We define a $4 \times 8$ matrix $M$ with entries in
 ${\mathbb{F}}_{2}$  by 
$$ M = \left[
\begin{array}{cccccccc} 
1 & 1 & 1 & 1 & 0 & 0 & 0 & 0 \\
0 & 0 & 1 & 1 & 1 & 1 & 0 & 0 \\
0 & 0 & 0 & 0 & 1 & 1 & 1 & 1 \\
1 & 0 & 1 & 0 & 1 & 0 & 1 & 0 \\
\end{array}
\right].$$
Let $C\subset {\mathbb{F}}_{2}^{8}$ be the row space of $M$. 
It can be easily verified that the row vectors of $M$ are 
linearly independent. Hence dim($C$) = 4. Throughout the paper
we will denote this code by $C_{8}$.

\medskip
If ${\bf v},{\bf w}$ are two elements of $\V$, then we define their product
${\bf v}\times {\bf w}$ by
${\bf v}\times {\bf w} = (v_{1}w_{1},\ldots ,v_{d}w_{d})$.
It is easy to see that ${\bf v}\cdot {\bf w} = 0$ if and only if
$|{\bf v}\times {\bf w}|$ is an even integer. 
\begin{prop}\label{basic}
Let $C\subset \V$ be a binary linear code.
\begin{enumerate}
\item{If $C$ admits a self-orthogonal doubly even basis $A$,
then $C$ itself is self orthogonal and doubly even.}
\item{$\mbox{dim}(C) + \mbox{dim}(C^{\perp}) = d$.}
\end{enumerate}
\end{prop}
{\bf Proof.}
Let ${\cal A}$ be the collection of all self orthogonal
doubly even subsets of $\V$ which contains $A$, and 
let $W$ be a maximal element of ${\cal A}$. Let ${\bf w}, {\bf w}^{'}$ be any
two elements of $W$. Since the dot product is a bilinear 
form on $\V$, the set $W \cup \{{\bf w} + {\bf w}^{'}\}$ is self
orthogonal. We note that for any two elements
${\bf v},{\bf w}\in \V$,
$$ |{\bf v}+{\bf w}| = |{\bf v}| + |{\bf w}| - 2|{\bf v}\times {\bf w}|.$$
Hence $W \cup \{{\bf w} + {\bf w}^{'}\}$ is doubly even.
By the maximality of $W$, ${\bf w} + {\bf w}^{'}\in W$. Therefore
$W$ is a subspace of $\V$. Since $A$ is a basis of $C$, this proves 
the first part. The second assertion is a consequence of the non-degeneracy
of the dot product.
$\hfill \Box$.

\medskip
Now we introduce a notion of non-degeneracy on binary linear codes.
For any $d\ge 1$ we define a map
$B : {\Z}^{d}\times \V \rightarrow \Z$ by
$$ B({\bf n}, {\bf v}) = \sum_{i : v_{i} = 1}n_{i}.$$
{\it Definition :\/}
Let $C\subset \V$ be a binary linear code. Then $C$ is 
said to be {\it integrally non-degenerate \/} if for 
all non-zero ${\bf n}\in {\Z}^{d}$, there exists
a ${\bf v}\in C$ such that 
$B({\bf n}, {\bf v})$ is non-zero.
\begin{prop}\label{e8}
For every $d\ge 8$ there exist proper subspaces $C, C^{'}\subset \V$
such that $C$ is an integrally non-degenerate code containing ${\bf 1}$
and for any two ${\bf x},{\bf y}$ in $C$, their product 
${\bf x}\times{\bf y}$ lies in $C^{'}$.
\end{prop}
{\bf Proof :}
First we will consider the case when $d =8$. We claim that 
the pair $(C_{8},E_{8})$ has the required properties.
It is easy to see that ${\bf 1}\in C_{8}$. 
Let $A\subset {\mathbb{F}}_{2}^{8}$ denote the set of row vectors of
the matrix $M$, as defined in Example 2. It is easy to check that 
$A$ is a doubly even  self orthogonal basis of $C_{8}$.
 By the previous proposition, $C_{8}$ is doubly 
even and self orthogonal.  In particular for any two 
${\bf x},{\bf y}\in C_{8}$, ${\bf x}\cdot {\bf y} = 0$ i.e.
${\bf x}\times {\bf y}\in E_{8}$.

Let ${\bf n}$ be a non-zero vector in ${\Z}^{d}$. Clearly, we can 
choose $i,j$ in $\{ 1,\ldots ,8\}$ such that 
$n_{i} + n_{j} \ne 0$. Let $\phi$ be the vector space
homomorphism from $C_{8}$ to 
${\mathbb F}_{2}^{2}$ defined by
$\phi({\bf v}) = (v_{i},v_{j})$ and let $E\subset C_{8}$
be the set defined by
$$ E = \{ {\bf v}\in C_{8}\ |\ \phi({\bf v}) = (1,1)\}.$$
It is easy to see that for every $k\in \{ 1,\ldots ,8\}$
there exists a vector ${\bf v}^{'}$ in $A$ such that 
${\bf v}^{'}_{k} = 1$. Hence we can choose vectors ${\bf v},{\bf w}$ 
in $A$
such that ${\bf v}_{i} = 1$ and ${\bf w}_{j} = 1$. It is easy
to see that there 
exists ${\bf x}$ in the set $\{ {\bf v}, {\bf w}, {\bf v} + {\bf w}\}$ 
which lies in $E$. Since $C_{8}$ has
16 elements and ${\mathbb{F}}_{2}^{2}$ has 4 elements,
the set $E = {\bf x} + \mbox{ker}(\phi)$ contains at
least 4 elements. Since $C_{8}$ is doubly even this implies that
we can find 
two distinct vectors ${\bf v}, {\bf w}$ in 
in $E$ such that $|{\bf v}| = |{\bf w}| = 4$. As 
${\bf v}\cdot {\bf w} = 0$, this
shows that $({\bf v}\times {\bf w})_{k} = 1$
if $k = i,j$ and zero otherwise. In particular,
$B({\bf n},{\bf v}\times{\bf w}) = n_{i} + n_{j}\ne 0$.
We also note that the map $B$ satisfies the identity
$$ 2 B({\bf m},{\bf x}\times {\bf y}) = 
 B({\bf m},{\bf x}) +  B({\bf m},{\bf y})
-  B({\bf m},{\bf x} + {\bf y}).$$
Hence we can find a vector ${\bf v}_{0}$ in the set 
$\{ {\bf v}, {\bf w}, {\bf v} + {\bf w}\}\subset C_{8}$ such that
$B({\bf n},{\bf v}_{0})$ is non-zero.
This proves the claim.
For $d > 8$, we define $C , C^{'}\subset \V = 
{\mathbb{F}}_{2}^{8}\oplus {\mathbb{F}}_{2}^{d-8}$ by
$C = C_{8}\oplus {\mathbb{F}}_{2}^{d-8}$ and
$C^{'} = E_{8} \oplus {\mathbb{F}}_{2}^{d-8}$. It is easy to verify
that the pair $(C,C^{'})$ has the desired properties.
$\hfill \Box$ 

\section{Non rigid actions}
As in Section 2, for any $d\ge 1$ by $X_{d}$ we denote the group 
$({\Z/2\Z})^{{\Z}^{d}}$
and by $S$ we denote shift action of ${\Z}^{d}$ on $X_{d}$. 
For any $x, y\in X_{d}$ we define their product $x \star y \in X_{d}$
by $x \star y ({\bf i}) =  x({\bf i})\, y({\bf i})$. It is easy to see that
$X_{d}$ becomes a compact topological ring with respect to this product and
$S({\bf i})(x \star y) = S({\bf i})(x) \star  S({\bf i})(y)$ for all
${\bf i}$ in ${\Z}^{d}$.

\medskip
The following proposition is the basis of our construction.
\begin{prop}\label{key}
Let $H, K\subset X_{d}$ be  proper Markov subgroups 
such that the actions $(H,S)$ and $(K,S)$ are mixing,
and for all $x,y$ in $H$, $x\star y\in K$. We define a
${\Z}^{d}$-action $(X,\alpha)$ and a map $f : X\rightarrow X$  by 
$$ (X,\alpha) = (H,S)\times (H,S)\times (K,S),\ 
f(x,y,z) = (x,y, x\star y + z).$$ 
Then $(X,\alpha)$ is a mixing zero entropy action of ${\Z}^{d}$, and 
 the map $f$ is a non-affine homeomorphism which preserves the
Haar measure on $X$ and commutes with the action $\alpha$.
\end{prop}
{\bf Proof :}
Since $H$ and $K$ are proper subgroups of $X_{d}$, 
by Proposition \ref{mar} both $(H,S)$ and $(K,S)$ have zero entropy. 
Since both $(H,S)$ and $(K,S)$ are mixing by our assumption, it
follows that $(X,\alpha)$ is mixing and has zero entropy.
 It is easy to see
that $f$ is a homeomorphism which commutes with the action $\alpha$.
From the standard results on skew products if follows that $f$ preserves
the Haar measure on $X$. So it remains to show that $f$ is a 
non-affine map. Suppose this is not the case. Comparing the last
coordinate we see that there exists a constant $c_{0}\in K$ and 
homomorphisms $\theta_{1}, \theta_{2} : H\rightarrow K$ and 
$\theta_{3} : K\rightarrow K$  such that
$$ x\star y + z = c_{0} + \theta_{1}(x) + 
\theta_{2}(y) + \theta_{3}(z) \ \ \forall x,y,z.$$
Putting $x= y = 0$ we see that $c_{0} = 0$ and 
$\theta_{3} = \mbox{Id}$. Putting $x = 0$ 
(resp. $y = 0$) we see that $\theta_{2} = 0$ (resp. $\theta_{1} = 0$).
Hence $x\star y = 0$ for all $x$ and $y$. On the other hand,
$x\star x\ne 0$ for any non-zero $x$. This contradiction completes
the proof. $\hfill \Box$

\medskip
For any binary linear code $C\subset \V$ we define a Markov subgroup
$X_{C}\subset X_{d}$ by
$$X_{C} = \{ x\in X_{d}\ |\ 
(x({\bf i} + {\bf e_{1}}),\ldots , x({\bf i} + {\bf e_{d}}))\in C
\ \ \forall {\bf i}\in {\Z}^{d} \},$$
where ${\bf e}_{1},\ldots , {\bf e}_{d}$ are the standard unit 
vectors in ${\Z}^{d}$.
The ideal $I(X_{C})$, 
as defined in Section 2, can be
described as follows : 
For any ${\bf v}$ in $\V$ 
we define a polynomial $p_{{\bf v}}$ in
${\cal R}_{d}^{(2)}$ by
$$ p_{{\bf v}} = \sum_{j=1}^{d} v_{j} u_{j}.$$
We note that $p_{\bf v}\cdot x = 0$ for any
${\bf v}\in C^{\perp}$ and $x$ in $X_{C}$. 
Since $(C^{\perp})^{\perp} = C$ by Proposition \ref{basic},
it follows that the $I(X_{C})$ is the ideal generated by the 
set  $\{ p_{{\bf v}}\ |\ {\bf v}\in C^{\perp} \}$.
As $p_{{\bf v} + {\bf w}} =  p_{{\bf v}} +  p_{{\bf v}}$
for any ${\bf v},{\bf w}\in \V$, we see that for any 
basis $A$ of $C^{\perp}$ the set 
$\{ p_{{\bf v}}\ |\ {\bf v}\in A \}$ generates the ideal
$I(X_{C})$. 

\medskip
{\bf Examples :}
If $C = E_{d}$  then $I(X_{C})$ is the principal ideal generated by 
$u_{1} + \cdots + u_{d}$. Since $ C_{8} = C_{8}^{\perp}$ and
the row vectors of the matrix $M$  form
a basis of $C_{8}$, it follows  that $I(X_{C_{8}})$
is the ideal  $<p_{1}, p_{2},p_{3},p_{4}>$, where
$p_{1}, p_{2},p_{3},p_{4}$ are given by
$$\begin{array}{lll}
p_{1} & = & u_{1} + u_{2} + u_{3} + u_{4}\\
p_{2} & = & u_{3} + u_{4} + u_{5} + u_{6}\\
p_{3} & = & u_{5} + u_{6} + u_{7} + u_{8}\\
p_{4} & = & u_{1} + u_{3} + u_{5} + u_{7}\\
\end{array}$$
\begin{lemma}\label{prime}
Let $A\subset \V$ be a subset and let 
$I_{A}\subset {\cal R}_{d}^{(2)}$ be the ideal generated by the set
$\{ p_{{\bf v}}\ |\ {\bf v}\in A \}$. Then $I_{A}$ is a prime
ideal.
\end{lemma}
{\bf Proof.}
Let $\pol$ be the polynomial ring in $d$ variables with
coefficients in $F_{2}$. We can identify $\pol$ with the subring
of $\lar$ which consists of all $p$ such that 
$c_{p}({\bf n}) = 0$ whenever $n_{i} < 0$ for some $i$.
Clearly for each ${\bf v}$ in $\V$, $p_{{\bf v}}$ lies in
$\pol$. For any  set $B\subset \V$ let
 $I_{B}^{'}\subset \pol$ be the ideal generated 
by the set $\{ p_{{\bf v}}\ |\ {\bf v}\in B \}$. 
We claim that $I_{B}^{'}$ is a prime ideal of $\pol$.
To prove this we note that $I_{B}^{'} = I_{C}^{'}$, where
$C\subset \V$ is the subspace generated by $B$.
If $\mbox{dim}(C) = k$, we define a subspace 
$C_{1}\subset \V$ by
$$ C_{1} = \{ (v_{1},\ldots ,v_{d})\ |\ v_{i} = 0 \ \forall i > k\}.$$
Let $\theta$ be a linear automorphism of $\V$ such that 
$\theta(C) = C_{1}$ and let $\overline{\theta}$ be the
automorphism of $\pol$ satisfying $\overline{\theta}(p_{{\bf v}})
= p_{\theta({\bf v})}$ for all ${\bf v}$ in $\V$.
Then $\overline{\theta}(I_{C}) = I_{C_{1}}$. It is easy to see that
$\pol / I_{C_{1}}$ is isomorphic to $F_{2}[u_{k+1},\ldots , u_{d}]$.
Therefore $I_{C_{1}}$ is a prime ideal. 
Since $\overline{\theta}(I_{C}) = I_{C_{1}}$ and 
$\overline{\theta}$ is an automorphism of
$\pol$, this proves the claim.

For any $N > 0$ we define ${\Z}^{d}_{N}\subset {\Z}^{d}$ by
$$ {\Z}^{d}_{N} = \{ (n_{1},\ldots ,n_{d})\in {\Z}^{d}\ |
\ n_{i} \ge N \ \forall i\}.$$
We observe that for any element $p$ in $\lar$,
$p$ lies in $\pol$ if and only if there exists $N > 0$ such that
$u^{{\bf n}}p \in \pol$ for all  ${\bf n} \in {\Z}^{d}_{N}$. 
Similarly for any $p$ in $I_{A}$,
$p$ lies in $I_{A}^{'}$ if and only if 
there exists $N > 0$ such that
$u^{{\bf n}}p \in I_{A}^{'}$  for all ${\bf n} \in {\Z}^{d}_{N}$.
Let $p_{1}, p_{2}$ be two elements of $\lar$ such that
$p_{1}p_{2}\in I_{A}$. Then we can  choose ${\bf n}\in {\Z}^{d}$
such that $u^{{\bf n}}p_{1}, u^{{\bf n}}p_{2} \in \pol$
and $u^{2{\bf n}}p_{1}p_{2}\in I_{A}^{'}$. 
By the above claim, $I_{A}^{'}$ is a prime ideal of $\pol$.
Hence either $ u^{{\bf n}}p_{1}\in I_{A}^{'}$ or 
 $ u^{{\bf n}}p_{2}\in I_{A}^{'}$. This implies that 
either $ p_{1}\in I_{A}$ or 
 $p_{2}\in I_{A}$, which proves the given assertion.
$\hfill \Box$
\begin{lemma}\label{mix}
Let $C\subset \V$ be a binary linear code.
\begin{enumerate}
\item{ If $C$ is integrally non-degenerate and contains ${\bf 1}$
then the ${\Z}^{d}$-action $(X_{C},S)$ is mixing.}
\item{The action $(X_{C},S)$ has zero entropy if and only if
$C$ is a proper subspace of $\V$.}
\end{enumerate}
\end{lemma}
{\bf Proof :}
1) For any ${\bf v}\in \V$, 
let $\phi_{{\bf v}}$ be the unique homomorphism from
${\cal R}_{d}^{(2)}$ to ${\cal R}_{1}^{(2)} = 
{\mathbb{F}}_{2}[z, z^{-1}]$ such that
$\phi_{{\bf v}}(u_{i}) = z$ if $v_{i} = 1$ and 
$\phi_{{\bf v}}(u_{i}) = 1$ if $v_{i} = 0$.
We choose ${\bf v},{\bf w}\in \V$ such that
${\bf v}\in C^{\perp}$ and ${\bf w}\in C$.
Since $C$ contains ${\bf 1}$ and ${\bf v}\cdot {\bf w} = 0$, 
$|{\bf v}|$ and $|{\bf v}\times {\bf w}|$ are even integers.
 Hence the sets 
$\{i\ |\ v_{i} = w_{i} = 1\}$ and 
$\{i\ |\ v_{i} = 1, w_{i} = 0\}$ contain even number of elements.
Therefore
$$\phi_{{\bf w}}(p_{{\bf v}}) = \sum v_{i}\phi_{{\bf w}}(u_{i}) = 0. $$
This shows  that
$I({X_{C}})\subset \mbox{ker}(\phi_{{\bf w}})$ 
for all ${\bf w}$ in $C$.
Let ${\bf n}$ be any non-zero element of ${\Z}^{d}$. As $C$
is integrally non-degenerate, there exists a ${\bf w}$ in $C$
such that $B({\bf n},{\bf w})$ is non-zero.
Since $\phi_{{\bf w}}(u^{{\bf n}}) = z^{B({\bf n},{\bf w})}$, 
we conclude that
$u^{{\bf n}} - 1$ does not lie in $I({X_{C}})$. By 
the previous lemma and 
Proposition \ref{mar} the action $(X_{C},S)$ is mixing. This proves 1).
The second assertion is an immediate consequence of Proposition \ref{mar}.
$\hfill \Box$  

\medskip
Now we turn  to the proof of Theorem \ref{th}

\medskip
{\bf Proof of Theorem \ref{th} :}
Let $C, C^{'}$ be proper subspaces of $\V$ satisfying the conditions
stated in Proposition \ref{e8}. 
Since $C$ contains ${\bf 1}$ and  
${\bf x}\times {\bf y}\in C^{'}$
for all ${\bf x},{\bf y}$ in $C$, it follows that $C\subset C^{'}$.
Hence $C$ and $C^{'}$ are integrally non-degenerate codes containing 
${\bf 1}$. Since $C^{'}$ is a proper subspace of $\V$, by
Lemma \ref{mix} the actions $(X_{C},S)$ and $(X_{C^{'}},S)$
are mixing and have zero entropy. 
It is easy to verify that for any $x,y$ in $X_{C}$, $x\star y$ 
is an element of $X_{C^{'}}$.
Now Theorem \ref{th} follows from Proposition \ref{key}.
$\hfill \Box$

\bigskip
\noindent
Address : School of Mathematics, Tata Institute of Fundamental 
Research,
Mumbai 400005, India.

\bigskip
\noindent
e-mail : siddhart@math.tifr.res.in 

\end{document}